\documentclass[11pt,a4paper]{article}
\usepackage[latin1]{inputenc}
\usepackage[english,francais]{babel}
\usepackage{amssymb,amsthm,amsmath,xspace}
\usepackage[mathscr]{eucal}
\usepackage{amscd}
\usepackage{pstricks}
\usepackage{multicol}

\newfont{\df}{cmssbx10}
\newfont{\rk}{cmsl11}
\newfont{\pf}{cmsl11}

\date{\empty}
\newcommand{\dis}{\displaystyle}

\begin{document}

\vskip -2,5truecm \centerline {\LARGE   Rotation numbers in Thompson-Stein
  groups } 

\vskip 0 truecm \centerline {\LARGE and applications}
\bigskip

\medskip

\centerline {{\rk  Isabelle LIOUSSE}} 

\bigskip

\bigskip

\bigskip

{\bf Abstract.} {\footnotesize We study the properties of  rotation numbers 
 for some  groups of piecewise linear homeomorphisms of the circle.   We use
 these  properties to  obtain  results on   PL  rigidity,  non isomorphicity,
 non exoticity  of  automorphisms,   non smoothability  for   Thompson-Stein
 groups. }

 \medskip


{\bf Mathematical subjet classification :}  37E10, 37C15, 20E

 \section {Introduction.}

In   1965, R. Thompson discovered the first example of a  finitely presented
infinite simple group. This group called $T$ is defined as a group of PL
homeomorphisms of the circle.

\medskip

We recall that an homeomorphism $f$ of the circle $S_r = 
\frac{\mathbb{R}}{\mathbb{Z}}$  of length $r$  is {\df piecewise linear (PL)}
 if there exists a finite subdivision  : $0<a_1<a_2......<a_p=r$ such that 

${\tilde f}_{\vert [a_i,a_{i+1}] } (x) = \lambda_i x + \beta_i$, where ${\tilde f}$ is a lift of $f$ to $\mathbb{R}$.

\medskip

The $a_i$'s are called {\df break points (or breaks)} 
 of $f$, their set is denoted by
$BP(f)$.

\smallskip

The $\lambda_i$'s  are called {\df slopes of $f$,} their set is denoted by
$\Lambda (f)$.

\smallskip

The ratio $\sigma_f(a_i) = {\frac {\lambda_i}{\lambda_{i-1} }}$ is called the
{\df  jump} of $f$ at $a_i$, their set is denoted by $\dis \sigma(f)$.

\medskip

The set of all orientation preserving (i.e. slopes are positive) PL
homeomorphisms of $S_r$ is  denoted by $PL^+ (S_r)$.

\bigskip

\noindent The {\df Thompson group $T$}  consists of  all orientation preserving 
homeomorphisms of ${\frac  {\mathbb{R} } { \mathbb{Z}}}$  such that :

-  breaks and their images are dyadic rational numbers,

-  slopes are  powers of $2$.

\bigskip

The group $T$ and its subgroup $F$ --consisting of the elements of $T$ that
fix $0$-- occur in many areas of mathematics : logic, homotopy theory, group
theory, dynamical systems. The standard survey  for Thompson groups is  [8].  
Dynamical properties
of $T$ and $F$ have been studied by Ghys and  Sergiescu in [11] 
 and Brin in [3]. 

\medskip

In  particular in [11],  
  it is  proved that :

- any non trivial representation $\phi : T \rightarrow Diff^r _+(S^1)$, $r\geq
2$ is semi-conjugate to the standard representation in $ PL^+ (S^1)$ ;

- $T$ admits faithful representations in   $Diff^\infty _+(S^1)$ and both of
the two following situations can occur : all $\phi (T)$-orbits are dense 
 or  $\phi (T)$ has an exceptional minimal.

\smallskip

As corollary and  using Denjoy theorem, it's shown in [11] 
that the rotation number of any element of $T$ is rational.

\medskip

In [12],  
 \'E. Ghys asked about the following problem : consider a piecewise
linear  homeomorphism $f$ of $S^1$ such that all its slopes, breaks, images of
breaks  are rational numbers, is it true that the  rotation number of 
$f$ is rational ?

This question has a negative answer. Counterexamples can be found in [2], a
Boshernitzan paper and not explicitly stated by 
 Herman in [13] chapt VI section7.

\bigskip

In [3],  
 M. Brin proved that the outer automorphism  group of $T$ has  
order 2, as is the group $OUT (Homeo^+ (S^1))$. This fact confirms the fact that the
finitely generated group $T$ is imitated the behaviors of the continuous group 
$Homeo^+ (S^1)$ (simplicity, cohomology, automorphisms).

An other interesting property of $T$ and $F$ due to Farley  is that $T$ and
$F$ are a-T-menable (see [10] ). 
  It's still unknown if $F$ is amenable or not.

\medskip

Moreover, generalisations of Thompson groups have been quite a lot defined and
studied  : Higman [14], Bieri-Strebel [1],  Brown [7],  Stein [22] and 
  Brin-Guzman ([6]). In this paper, we  consider the
  generalisations in  $PL(S^1)$ and follow the descriptions  
of Bieri-Strebel and Stein.

\bigskip

\bigskip

\noindent {\bf  Thompson-Stein groups.}

\medskip

{\bf Definitions.} Let $\Lambda \subset \mathbb{R}^{+*}$ be a multiplicative
subgroup,   and  $A \subset \mathbb{R}$ be an additive subgroup, invariant by
multiplication by   elements of $\Lambda $.

\noindent We define $T_{r,\Lambda,A}$ as the subgroup of  $PL^+ (S_r)$
consisting of elements with slopes in $\Lambda$, breaks and their images in
$A$. The slope set of  $T_{r,\Lambda,A}$ is $\Lambda$ and  the break 
set of  $T_{r,\Lambda,A}$ is $A$.

\medskip

 In the case  $\Lambda=<n_i>$ the multiplicative group generated by $p$
  integers $1<n_1<...<n_p$  independent (in the sense of free abelian group
  generators, this  is equivalent to the statement that the $\log n_i$ are
  $\mathbb{Q}$-independent) and $A = \mathbb{Z}[\frac{1} { n_1},...,\frac{1} {n_p}]=
  \mathbb{ Z }[\frac{1}{m}]$ where $m=lcm (n_i)$ the least common multiple of the $n_i$ : 

\smallskip

the group  $T_{r,\Lambda,A}$ is denoted by $T_{r, (n_i)}$ and called a {\df
Thompson-Stein group}. It's the smallest non trivial  subgroup
$T_{r,<n_i>,A}$.  The  subgroup of $T_{r, (n_i)}$ consisting of  elements that
fix $0$  is denoted  by $F_{r, (n_i)}$.

\medskip

When $p=1$, the group $T_{r,(m)}$ [resp. $F_{r,(m)}$]  is denoted by $T_{r,m}$
[resp. $F_{r,m}$]. The classical Thompson groups are $T=T_{1,2}$ and $F= F_{1,2}$.

\bigskip

\noindent {\bf Questions.}  For these groups many natural questions arise :

\noindent How are the rotation numbers of the  elements of $T_{r, (n_i)}$ ?

\noindent How are automorphisms of $T_{r, (n_i)}$ ?

\noindent  Is it possible to describe all 
$T_{r, (n_i)}$ up to isomorphisms, up to quasi isometry ?

\noindent Is it possible to classify the representations of  $T_{r, (n_i)}$
in $PL^+ (S_r)$, $Diff ^k _+ (S_r)$ ?   In particular,  are the  $T_{r, (n_i)}$
(infinitely) smoothable ?

\bigskip

{\bf Properties.}  M. Stein in  [22] proved that :

\noindent- $T_{r, (n_i)}$ and    $F_{r, (n_i)}$  are finitely presented and $FP^\infty$,

\noindent-  $T''_{r, (n_i)}$ and    $F'_{r, (n_i)}$ are simple groups,

\noindent - precise homology informations.

\medskip

In  [1] --the author is attached  to mention that she
couldn't get  this unpublished paper, so she refered to results  that are
cited or reproved in [22], [3] and [6]--  Bieiri and Strebel 
have studied isomorphisms and automorphisms of $T_{r,\Lambda, A}$.

\smallskip

They've described PL isomorphisms between distincts $T_{r, (n_i)}$ ; in
particular,  they
proved that : ``  $T_{r, (n_i)}$  and   $T_{r', (n_i)}$  are  PL-isomorphic
if and only if , $(r-r')$ is an integer multiple of $d= gcd(n_i-1)$  the highest common divisor of
the $n_i-1$''.

\smallskip

They also proved that automorphisms of $F_{r, \Lambda, A}$  are not 
 exotic (i.e. they  are  realized as PL conjugation) provided 
 $\Lambda$ is a dense
multiplivative group. 

\smallskip

To extend this result  to $T_{r, (n_i)}$ (in particular for $p\geq 2$), it's reasonnable to think about  localization results (such
results can be found in  [3] and [6] for $T_{r,m}$ when 
$r$ is  not a multiple  of $d=m-1$). A very surprising theorem 
is proved in [6] : $T_{m-1, m}$ has (infintely many) exotic automorphisms.

\medskip

Here, we developp a  different approach  to study  the
automorphisms of Thompson-Stein groups : 
 we analyse the rotation numbers and  use  a rigidity result
for actions in $PL^+ (S^1)$ or $Diff^k_+(S^1)$ to  obtain
as corollaries  informations on  isomorphisms and automorphisms.

\eject

\noindent{\bf \large Theorem 1. Rational rotation numbers.} {\it Let    
$m\geq 2$, $r\geq 1$ and    $q\geq 1$ be integers. Then}

\medskip

{\it A. Each    element of $ T_{r,<m>,\mathbb{Q}}$  has rational  rotation
  number, so has periodic points.}

\medskip

{\it B. The following   properties are equivalent :

{\df (1)  there exists $p\in \mathbb{N}^*$ such that $ T_{r,m}$
  contains  an  element of reduced rotation number  $\frac{p}{q }$,}

(1') for all  $p\in \mathbb{N}^*$,  the group $ T_{r,m}$
  contains  an  element of reduced rotation number $ \frac{p}{q }$,

{\df  (2) gcd $(m-1, q)$ is a divisor of $r$,}

(3) the group $ T_{r,m}$
  contains an element of order $ q$ 

(3') for all  $p\in \mathbb{N}^*$, the group $ T_{r,m}$
  contains an element of order $ q$ with  rotation number  $ \frac{p}{q}$.}

\bigskip

\bigskip

\noindent {\bf \large Immediate consequences of Theorem 1.} Fix  two integers
$m\geq 2$ and  $r\geq 1$ .

\noindent {\df  On the rotation number 
map } $ \rho : {\text { Homeo}} ^+ (S_r) \rightarrow \frac {\mathbb{R}}{\mathbb
{Z}}$ defined by   $\dis  \rho(f) = \lim\limits_{n\rightarrow \infty} ({\tilde
f^n (0)/ rn }  ) \ (mod 1)$.

\smallskip

{\it 1. The image  set $\rho ( T_{r,m}) \subset \frac{\mathbb{Q}} {\mathbb {Z}}$ and any 
  number  in  $\rho (T_{r,m})$  is  realized as 
the rotation number of a finite order element of $ T_{r,m}$. 

\smallskip

 2.  Any rational number is realized as 
the rotation number of an element of finite order  of $ T_{m-1,m}$ ; in
particular,  $\rho ( T_{m-1,m}) =  \frac{\mathbb {Q} }{\mathbb {Z}}$.

\smallskip

3. The group $ T_{m-1,m}$  is the only group 
for $0< r <m-1$ that contains elements of order $m-1$. 
In particular,  $\rho ( T_{r,m}) \not=  \frac{\mathbb{Q}}{\mathbb{Z}}$  for $0< r
<m-1 \ (mod \  m-1)$.

\smallskip

 4. If $m,r$ are  odd, the group $ T_{r,m}$ doesn't
  contain elements of order   2 (of any even order).

\smallskip

 5. The image set $\rho (T_{1,m}) = \{ \frac{p}{q} \ : \gcd (q, m-1)=1\}$. If $m-1$ is a prime number then $\rho (T _{ 1,m }) =\rho (T _{ r,m })$, 
for all $0< r <m-1 \ (mod \  m-1)$.

\medskip

\noindent{\df On  the isomorphic problem.}

{6. The group  $ T_{m-1,m}$  is  isomorphic to none of the groups 
 $ T_{r,m}$ with $0\leq r <m-1$.}

 7. If $m$ is odd,   $r$ odd  and  $r'$ even,  the groups 
 $ T_{r,m}$  and  $ T_{r',m}$ are not isomorphic.}
\bigskip

\noindent In the next  theorem, we prove that  Thompson-Stein groups contain 
free abelian groups that act freely on $S_r$. The  rank of such a group 
 depends on the rank  of  the  slope set. 

\medskip

\medskip

\noindent{\bf \large Theorem 2. Irrational rotation numbers.}
 {\it  Let $G= T_{r,(n_i)}$ be a Thompson-Stein group associated  to $p\geq 1$
  independent integers $n_i$. Then : }

\medskip

{\it A. The group  $G$ contains  a free abelian group of rank $p-1$  acting
  freely. That is equivalent to the statement that  $G$ contains   $(p-1)$
  commuting homeomorphisms with    irrational  rotation numbers $\rho_i$ such that  $1$ and the $\rho_i$ are  $\mathbb{Q}$-independent.}

\smallskip

{\it B.  The group $G$ doesn't contain  any  free abelian group of rank $p$
  acting  freely on the circle.}

\medskip

{\it C. The rotation number of any element 
  in  a free abelian subgroup  of rank $2$ of $G$  is written in the form 
 $\rho(g) = \frac{\log \alpha } { \log \beta}$, with $\alpha, \beta \in <n_i>.$}

{\it C'.   Any number $\frac{\log \alpha }{ \log \beta}\in [0,1[$  with
  $\alpha,  \beta \in <n_i>$ can be realized as the rotation number of an
  element of   $ T_{d,(n_i)}$, where  $d= gcd (n_i-1)$.}

\medskip

\medskip

\noindent {\bf \large Consequences of theorem 2.} {\it Let $G= T_{r,(n_i)}$
  and  $A= \mathbb{Z}[\frac{1}{ m}]$, with $m=lcm(n_i)$.}

\smallskip

{\it 1. The  free rank of $G$ defined   by  $frk(G)=
\max \{n \in \mathbb{N} : \exists  \mathbb{Z}^n \subset G$ acting freely$\} = p-1$.

\smallskip

   2. Any representation $\phi$ from $G= T_{r,(n_i)}$, $p\geq 3$  into $PL^+(S^1)$ 
that is topologically conjugate to  the standard representation in
$PL^+(S_r)$  is $PL$-conjugate to the standard representation.

\smallskip

3. The group  $G= T_{r,(n_i)}$, $p\geq 3$   has no exotic automorphism ; more
precisly, its automorphisms are in realized by conjugation by maps in 
$T_{r, inv A, A}$, where $inv A = \{a \in A : \frac{1}{a}\in A \}$. In
particular, $OUT G$  has order 2 provided the $n_i$ are coprime.

\smallskip

4. Let $r\in \mathbb{N}^*$, two distincts Thompson-Stein groups
$T_{r,(n_1,...n_p)}$ with $p\geq 3$  and $T_{r,(n'_1,...,n'_q)}$  with $q\geq 1$  are not isomorphic.}

\medskip

\noindent {\bf Remarks on theorem 2 and aknowlegments.} According to G. Rhin 
([21]), the  $\frac{\log \alpha}{\log \beta}$   are diophantine numbers. I'm grateful to Yann Bugeaud
and Nicolas Brisebarre  for indicating  me      references and 
 details on  this property.

\bigskip

The proof of theorem 2A is based on extensions  of Boshernitzan examples. Part B and C are    consequences   of   a rigidity result for
free PL-actions of  $\mathbb{Z}^2$ on the circle due to Minakawa  ([18]) and
 the explicit calculation of rotation  numbers  for some particular PL-homeomorphisms.

\bigskip

\noindent {\bf \large Theorem 3 (Representations in $PL(S^1)$ or $Diff(S^1)$). }

\noindent {\it Let $G= T_{1,(n_i)}$ be the Thompson-Stein group 
associated  to $p\geq 2$
  independent integers : $(2,n_2,...,n_p)$. Then : }

\smallskip

{\it A. Each  non  trivial [faithfull]  representation $\phi$ from $G$ into
  $PL(S^1)$  or  $Diff^2(S^1)$ is  topologically conjugate to the standard representation in $PL(S_r)$.

\medskip
 B. There exists  $k\geq 2$  depending on the $log n_i$'s diophantine
 coefficients  such that any  representation  from $G$ into
 $Diff^k (S^1)$ is trivial [has finite image]. In particular, $G$ is not
 realizable in  $Diff^\infty  (S^1)$.}

\bigskip

The proof  of theorem 3  is the combination of  : 

- dynamics properties  of $G$ and its  subgroups  acting on $S^1$, following  Ghys-Sergiescu
approach. It's important to say that the theorem K in [11] holds  for
PL-actions instead of $Diff^2(S^1)$ (as Koppel lemma holds  for PL-actions)
and  doesn't extend to $T_{r,m}$ with $m>2$  since the Euler class triviality
argument  doesn't work with $rk (H^2(T_{r,m}, \mathbb {Z})>2$.

- and the Brin version ([4]) of a Rubin's  theorem ([21]).


\bigskip

\noindent {\bf Consequence of theorems 2 and 3.} {\it  Let $G= T_{1,(n_i)}$ be a Thompson-Stein group associated  to $p\geq 3$
  independant integers : $(2,n_2,...,n_p)$. Any non  trivial representation $\phi$ from $G$ into $PL(S^1)$ is  PL-conjugate to the standard representation.}

\bigskip

\noindent A important  piece  for our  proofs is the  Bieri-Strebel criterion.

\medskip

{\bf Definition.} Two real intervals $[a,c] $ and $[a',c']$ are {\df $PL_{\Lambda,A}$-equivalent}  if
there  exists  a PL homeomorphism with slopes in $\Lambda$ and breaks in $A$
that takes $[a,c] $ to  $[a',c']$.

\medskip

{\bf   Bieri-Strebel criterion.}   The intervals  $[a,c] $ and $[a',c']$ are
$PL_{\Lambda,A}$-equivalent if and only if their length are equal modulo
$(1-\Lambda)A$, where the set  $(1-\Lambda)A = \{ \sum (1- \lambda_i) a_i,
\lambda_i \in \Lambda, a_i\in A\}$ (for a proof, see [22] appendix).

\medskip

{\bf   For Thompson-Stein groups.}  Suppose that   $\Lambda= <n_i>$  and
$A=\mathbb{Z}[\frac{1}{ m}]$ with  $m=lcm(n_i)$, then  :

-  the   set  $(1-\Lambda)A = dA= (d\mathbb {Z}) \Lambda$,  where  $d =gcd  (n_i-1)$.

-  two integers   are equal  modulo $(1-\Lambda)A$ iff 
their difference is a multiple of $d$.

-  the groups $T_{r,(n_i)}$ and  $T_{r',(n_i)}$ are PL-conjugate 
 provided $r=r'$ (mod d).

\bigskip

{\bf  Proof of   $(1-\Lambda)A = dA$. }

\smallskip

{\df Step 1.}  $(1-\Lambda)A \subset  dA$. As  $(1-\Lambda)A $ is a submodule, 
it suffices to prove that $(1-\lambda) a \in dA$, with $a\in A$ and $\lambda =
n_1^{s_1}  ...n_p ^{s_p} \in \Lambda$. By reducing to same denominator, $(1-\lambda) a = (1- n_1^{a_1}  ...n_p ^{a_p} ) a = 
( n_1^{b_1}  ...n_p ^{b_p}- n_1^{c_1}  ...n_p ^{c_p} ) a'$ with $b_i, c_i $
positive integers, $a'\in A$. By replacing the $n_i$ by $k_i d +1$ and
developping, we get    $(1-\lambda) a =  (dN) a' $ with  $N\in \mathbb{N}$.

\medskip

{\df Step 2.}  $dA\subset (1-\Lambda)A $. Let $a\in A$, using  Bezout
identity  we get  $d = u_1(n_1-1) + ....+ u_p(n_p-1)$ and hence  $da  =
\sum  (n_i-1)   (u_i a)  \in (1-\Lambda) A$.

\bigskip

\eject

\section{Proof of Theorem 1.}

\noindent {\bf\large  Proof of  theorem 1.A.} The  first step is to prove
that each $f \in T_{r,m}$ has rational  rotation number. The  homeomorphism $f\in T_{r,m}$  of  $S_r =\frac{[0,r]}{0=r}$
is identified  with  the bijection  $\tilde f$ (mod r) of  $[0,r[$, where
$\tilde f$ is a lift to $\mathbb{R}$ of $f$.

For all  $x \in [0,r[$, we have   
$\displaystyle f(x) = m^{e(x)} x +\frac{ p(x) } {m^{k(x)}} $, with  
  ${e(x)}\in \mathbb{Z}$, $ p(x)\in \mathbb{Z}$ and ${k(x)}\in \mathbb{N}$. Note  that
  these three functions   are  bounded by  constants independent  of $x$.

The orbit of  $0$  contains only  $m$-adic numbers, that is 
for  all $n\in \mathbb{Z}$, the point  $f^n(0)$ is written is the form 
\ \ \  $\dis f^n(0) = {\frac{M_n}{m^{N_n}}},$ \ \  with   $N_n \in \mathbb{N}$ and
$(*)$ \  $M_n\in \mathbb{N}$ is not a multiple of
$m$ or is zero. Remark that if  $m$ is not a prime number,  the fraction may  be
unreduced, but the  condition $(*)$ assumes  the uniqueness  of the expression
$f^n(0) = {\frac{M_n}{m^{N_n}}}$.

\medskip

If there  exists a non zero integer   $n$  such that 
  $M_n=0$ then  $f^n(0)=0$, therefore  the orbit of $0$  is  periodic and the rotation number of   $f$
is rational, we are done. 
 In what follows, we assume that   for all positive integer $n$ 
 the integer  $M_n$ is not a multiple of  $m$.   

\bigskip

\noindent {\bf Fact 1 :} {\it $N_n \rightarrow +\infty$ when $n
  \rightarrow +\infty$ or the orbit of 
$0$ is  periodic.}

\smallskip

By absurd, if the sequence  $N_n$ doesn't tend to  $+\infty$,  there  exists
a  bounded subsequence $N_{s_n}$, we write $N_{s_n}\leq B$. In this case, the
subsequence  $f^{s_n} (0)$  is  contained in a  finite set :

\noindent  $ \dis \{\frac {\displaystyle p } { \displaystyle  m^k }\ \ {\text { with }} \   
0\leq k \leq B \ {\text { and }} \  0\leq p\leq r.m^B \}$. But,  this is  possible iff the orbit of  $0$ is periodic.

\bigskip

\noindent  {\bf Fact 2 :} {\it If  $N_n \rightarrow +\infty$   then  $Df^n(0)
  \rightarrow 0$  when $n   \rightarrow +\infty$.}

\medskip

We suppose that  $N_n \rightarrow +\infty$  and compute   $\dis f^{n+1} (0)$ :

\smallskip

\noindent  $\dis f^{n+1} (0) = f(f^n(0)) 
= m^{e(f^n(0))} \frac{M_n } { m^{N_n}}+\frac{ p(f^n(0))} {m^{k(f^n(0))}} 
 =  \frac{M_n }{ m^{N_n - e(f^n(0))}}+ \frac{ p(f^n(0))} {m^{k(f^n(0))}}.$

\smallskip

\noindent As  $N_n $ tends to   $+\infty$,  as  $e$ and  $k$ are 
bounded, there is an integer  $n_0$ such that for all  $n\geq n_0$,  we have  
${N_n - e(f^n(0))}> k(f^n(0))$. 
It follows that   :

$$f^{n+1} (0) =
\frac  {\displaystyle  M_n + m^{ N_n - e(f^n(0)) -k(f^n(0)) }  p(f^n(0)) } {\displaystyle  m^{N_n - e(f^n(0))} }
=\frac{ \displaystyle M_{n+1} }{ \displaystyle  m^{N_{n+1}} }.$$

\smallskip

\noindent

As $m$ is  a divisor of   $m^ {N_n - e(f^n(0))- k(f^n(0))}$ but  not a divisor
of $M_n$, the integer  $m$  is not a divisor of 
$M_n+ m^ {N_n - e(f^n(0))- k(f^n(0))}p(f^n(0))$. The uniqueness of
the expression   $f^{n+1} (0)=\frac{ M_{n+1} }{m^{N_{n+1}}}$, implies that 
 $N_{n+1} = {N_n - e(f^n(0))}$, for all $n\geq n_0$.

\medskip

\noindent Thus, for all  $p\in \mathbb{N}^{*}$, we have  $\displaystyle 
N_{n_0+p} =   N_{n_0} - \left(e(f ^{n_0}(0))+...+e(f^{n_0+p-1}(0)) \right).$
     Hence, $\displaystyle e(f
^{n_0}(0))+...+e(f^{n_0+p-1}(0))
  \rightarrow -\infty$  when 
$p\rightarrow +\infty$.
\noindent Finally,   $\displaystyle e(0)+e(f(0))+...+e(f^{n}(0))  \rightarrow
-\infty$  when 
$n\rightarrow +\infty$. But  $\displaystyle Df^n(0) = m^{\displaystyle e(0)+e(f(0))+..
.+e(f^{n}(0)) }$ so  $\displaystyle Df^n(0) \longrightarrow 0$.

\bigskip

\noindent {\bf Conclusion.}  If the 
  rotation number of  $f$ is   irrational : 

\smallskip 

- The  $f$-orbit of  $0$ is not  periodic. From this and 
 fact   $1$  follows  that  $N_n\rightarrow +\infty$  and thus 
using fact 2 we get that $ Df^{n} (0)\rightarrow 0$.
\smallskip

- Denjoy's inequalities apply to $f$ and  give  $ Df^{q_n} (0)\geq
e^{-Varlog Df}$ that can't tend to  $0$.

\smallskip

\noindent From this   contradiction follows the rationality of 
the  rotation number of  $f\in T_{r,m}$. 
 
\bigskip

\noindent{\bf End of the proof of  theorem 1.A.} The group $ T_{r, <m>,\mathbb
  {Q}}$ is identified with  a  group of   
piecewise linear  bijections of $[0,r[$ with slopes 
power of $m$,  breaks and images of  breaks 
 rational numbers. Let $f \in   T_{r, <m>,\mathbb{Q}}$, on each continuity
 interval  of  $Df$, we write $f(x) = m^{e(x)}.x + \frac{p(x) } { q(x)}$. 
Let  $Q$ be the  lcm of the  $q(x)$, the homothety ${\cal H}_Q : 
S_r  \rightarrow S_{Qr}$ (${\cal H}_Q (x)=Qx$) conjugate $f$ to an element of  $T_{Qr, m}$
for which the first part applies.

\bigskip

\bigskip

\noindent {\large \bf Proof of theorem 1.B.}

\medskip

\noindent {\bf Easy  implications.}

$(3)\Rightarrow (1)$ is clear, since  $f$ of order $q$ has rotation
number $\rho(f) =\frac {p} { q}$ with $(p,q)= 1$.

$(1) \Leftrightarrow (1')$ and $(3) \Leftrightarrow (3')$  are obtained by
choosing suitable powers of $f$.

\bigskip

\noindent  {\bf We now prove the  two implications $(1)\Rightarrow
  (2)\Rightarrow (3)$.} 

\smallskip

{$(1)\Rightarrow (2)$.} By hypothesis, there exists $f\in T_{r,m}$ with
$\rho{(f)} =\frac{p} {q}$.  According to  Poincaré, $f$ has a periodic point $a$
of period $q$.  Replacing $f$ by a suitable power, we may assume  $\rho{(f)}
=\frac{1} { q}$.

Fix $\tilde f$ a lift of $f$ to $\mathbb{R}$ and identify  $f$  with  the
bijection $\tilde f (mod \  r)$ of $[0,r[$.  By replacing  $a$ by
 its nearest from $0$ iterate, the orbit of $a$ can be ordered as  :
$$0\leq a=a_0 <f(a)=a_1 <....< f^{q-1}(a)=a_{q-1}<r$$

\medskip

The intervals $[0,a]$ and $[f(0),f(a)]$ are 
$PL_{\Lambda,A}$-equivalent so, using the Bieri-Strebel criterion,  $a = (f(a) -f(0)) \  mod  (1-\Lambda)A$ 
 and thus  $ f(a) - a = f(0)  \  mod (1-\Lambda)A$.

\medskip

Furthermore, (via the map $\tilde f$) the intervals  $[a_0,a_1]$, $[a_1,a_2]$,
... and $[a_{q-1},a_0+r]$ are  $PL_{\Lambda,A}$-equivalent so their lengths $l_i$ 
   all  equal $ mod  (1-\Lambda)A$  the length of  the interval $[a_0,a_1]$
   that is \ $f(a) - a$.

Finally,  the lengths  $l_1$, $ l_2$,....$l_{q}$ of  the intervals
$[a_0,a_1]$, $[a_1,a_2]$, ...,$[a_{q-1},a_0+r]$  equal  $f(0)$ \ \   $mod \
(1-\Lambda)A=dA$. Adding $l_1+ l_2+....+l_{q}$,  we obtain  $r = qf(0) \  mod
\ dA$,  so $r - qf(0)   \in     (m-1) \mathbb{Z} [\frac{1} { m}]$. Thus,
there exists integers $u,v,s$ such that  $r - qf(0) =(m-1)\frac {v} { m^s}$ \ 
and  \ $f(0)  = \frac {u} { m^s}$, so such  that 
$ m^s r -qu = (m-1)v$, that is to say $ m^s r = qu + (m-1)v$.  This implies  that  $ m^s r$ is a multiple of  
$gcd(q, m-1)$.  But $(m-1)$ and $m$ (also $m^s$)  are  coprime,   we 
conclude therefore that   $ r$ is a multiple of   $gcd(q, m-1)$.

\bigskip

{$(2)\Rightarrow (3)$}.  Fix  $r$, $m$  two  positive integers  and suppose
that   $ r$ is a multiple of   $gcd(q, m-1)$. Using Bezout,
$r=uq + v(m-1)$, so $r= uq $ modulo $(m-1)$. The Bieri-Strebel  criterion implies that the groups $T_{uq, m }$ and $T_{r, m}$
are isomorphic. Furthermore,  the  group
$T_{uq, m }$ contains the   rotation $ x\mapsto x+u$ of order $q$. This
completes the proof.

\section{Proof of theorem 2 and consequences.}
The proof of theorem 2.A is based on generalisations of the examples
gived by Boshernitzan in [2]. The first type of generalisation is to 
pass from the usual circle to  $S_r$ and the second is to construct carefully
examples  in order to obtain  break points in the suitable set $A$ and 
commuting examples.

\subsection {Boshernitzan examples on  $S_r$.}

Consider  $S_r= \frac{\mathbb{R}} { r \mathbb{Z}}=\frac{[0,r]} {0=r}$ the circle of length  $r\in \mathbb{N}^*$ and $f \in PL ^+ (S_r)$ with exactly two break points  $a$ and
$f(a)=0$. The
map  $f$  is identified with  the  bijection $\tilde f$ (mod r)
 of  $[0,r[$.

\vskip 3 truecm
\hskip 6 truecm
\begin{picture}(0,0)(0,0)
\put(0,0){\line(1,0){80}}
\put(0,0){\line(0,1){80}}
\put(0,80){\line( 1,0){80}}
\put(80,0){\line(0,1){80}}
\put (0,40){\line(1,0){80}}
\put (0,40){\line(3,4){30}}
\put (30,0){\line(5,4){50}}
\put(30,0){\line(0,1){80}}

\put(-2,-2){$\bullet$}
\put(-10,-10){$\scriptstyle 0$}

\put(28,-10){$\scriptstyle a$}
\put(28,-2){$\bullet$}

\put(78,-10){$\scriptstyle r$}
\put(78,-2){$\bullet$}

\put(-20,38 ){$\scriptstyle f(0)$}
\put(-2,38){$\bullet$}

\put(12,68){$\scriptstyle \lambda_1$}

\put(53,27){$\scriptstyle \lambda_2$}

\put(-10, 78){$\scriptstyle r$}
\put(-2,78){$\bullet$}
\end{picture}
\bigskip

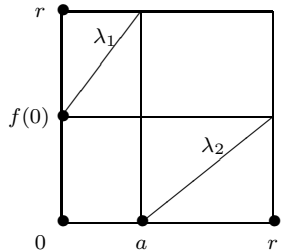
\centerline {\it Fig.  1}

\bigskip

 Fixing the initial break to be $0$, such a map is uniquely defined by its
slopes $\lambda_1$ and $\lambda_2 $.  Then, it admits $0$ , 
 $a = r \frac {1-\lambda_2} { \lambda_1-\lambda_2}$ as breaks and $0$ , $b= f(0) = 1- \lambda_1 a$ as 
images of breaks.

\medskip

Conversely, given two positive numbers $\lambda_1$ and $\lambda_2 $ such that 
$(\log  \lambda_1\log \lambda_1) <0$, there exists a unique $ f \in PL ^+ (S_r)$
with  breaks $0$ and $a$ and slopes  $\lambda_1$ and $\lambda_2 $ such that $ 
f(a)= 0$. This PL homeomorphism  $f$ is denoted by $f _{r, \lambda_1,  \lambda_2}$.   Let's define : 

\medskip

-  the homothety  ${\cal H}_r : \left\{ \begin{array}{ll} S^1=S_1  &
 \rightarrow S_r\\ \ \ \ \ x & \mapsto rx \end{array} \right. $ of ratio $r$ \ \ \    and 

- $ h _{\sigma}$  the  homeomorphism of $S^1$ identified with  the restriction to
 $[0,1[$ of one of its lift :

\ \ \ \ \ \ \ \ \ \ \ \ $\dis h_{\sigma} (x) =\frac {\sigma^x-1 } {\sigma -1 }$
\ if $\sigma\in
]0,1[$  \ \ \ \  and \ \ \ \   $ h_{1} (x)=x$. 

\bigskip

\noindent It's easy to check  (by an explicit calculation,
 as it's done in [2]), that  : 

\smallskip

$\triangleright$   the  map   ${\cal H}_r  \circ h_{\frac{\lambda_1} {\lambda_2} }$ 
 conjugates $f_{r, \lambda_1,  \lambda_2}$ to the rotation of $ S^1$ of angle
 $\frac   {\log  \lambda_1}  {\log \lambda_1-\log \lambda_2}$,

\smallskip

$\triangleright$  in particular, $f$  has   rotation number 
 $\dis \rho (f)=  { \frac{\log  \lambda_1} {\log \lambda_1-\log \lambda_2}}
 (mod \ 1)$.

\medskip

$\triangleright$ Conversely, for all $\rho \in S^1$, for all ${\sigma}\in
]0,1[$,   the map $({\cal H}_r  \circ h_{\sigma })\circ R_\rho \circ ({\cal
  H}_r  \circ h_{\sigma })^{-1}$ is in $PL^+(S_r)$ and has two breaks $0$ and
$a$ satisfying $f(a)=0$, it's called  a  {\df Boshernitzan on $S_r$.}

\bigskip

The map defined by $B_\rho^\sigma  =      h_{\sigma } \circ R_\rho \circ  h_{\sigma
}^{-1}$ for  ${\sigma}\in ]0,1]$  is  a  classical  Boshernitzan on $S^1$.

\bigskip

{\underline {To sum up}},\   \ $\dis f _{r, \lambda_1,  \lambda_2}=  ({\cal H}_r  \circ
h_{\sigma })\circ R_\rho \circ ( {\cal
  H}_r  \circ h_{\sigma })^{-1} = {\cal H}_r \circ  B _\rho ^\sigma\circ {\cal
  H}_r ^{-1},$

\smallskip

\noindent where $(\lambda_1,  \lambda_2)$  and $(\rho,  \sigma)$ are
related by  $\left \{\begin{array}{ll}& \sigma= \frac{\lambda_1} {\lambda_2}\\ 
& \rho ={\frac{\log \lambda_1} {\log\lambda_1- \log\lambda_2}}
\end{array}\right.$ \ \ \ \ and \ \ \
$\left \{ \begin{array}{ll} & \lambda_1 = \sigma ^{\rho -1}\\ 
& \lambda_1 =\sigma ^{\rho} \end{array}\right.$

\medskip

It's easy to check   that  two Boshernitzan of irrational rotation
numbers  commute if and only if they have  
the same  jump $\sigma=\frac{\lambda_1} {   \lambda_2}$ at $0$, 
since they are conjugate to rotations through  the same normalized homeomorphism.

\subsection {The Ghys question and remarks.}

Let's $n_1$ and  $n_2$  be two  independents integers.
Consider  $f$ the   Boshernitzan of  $S^1$ with slopes 
  $\lambda_1= n_1$ and  $ \lambda_2=\frac{1} {n_2}$. The data of $f$  are 
rational. Thus   $f\in T_{1,\mathbb{Q},\mathbb{Q}}$ but its  rotation number 
 $\rho (f) =\frac{\log n_1} { \log n_1+\log n_2}
\notin \mathbb{Q}$.   In general, because of the denominator in the formula that gives the break
point $a= \frac{1 - \lambda_2 }{\lambda_1 -\lambda_2}$,  the previous PL
homeomorphism is  not in the Thompson-Stein group
 $T_{1,(n_1,n_2)}$.   For example, by taking  $n_1= 2$, $n_2= 3$, the associated
 break point  is $a=\frac{2} { 5}\notin \mathbb{Z}[\frac{1}{6}]$. Note that, for proving that a   Boshernitzan $f$   of  $S_r$ is
in 
$T_{r, \Lambda, A}$, it  suffices to prove that its slopes 
$\lambda_1$ and  $\lambda_2$ are in  $\Lambda$ and that the break  $ a$ is in
$A$ (from this follows that $f(0)\in A$ because
$f(0)=r-\lambda_1a = \lambda_2 (r-a)$). The  examples  with irrational
rotation numbers  that we will construct  in Thompson-Stein groups are not 
Boshernitzan examples they are PL conjugate to  Boshernitzan examples.

\subsection{ Proof of theorem 2.A.}

Here, we construct elements of irrational rotation numbers in any 
Thompson-Stein group having slope set of rank  $p$ at least 2. Let $G =
T_{r,(n_1,...,n_p)}$ be a Thompson-Stein group, $m=lcm(n_i)$ and $d  =gcd(n_i-1)$.

\medskip

\noindent {\bf Step 1 : there exists $\Pi \in \Lambda$, $\Pi = n_1 ^{\alpha_1} n_2
^{\alpha_2}... n_p ^{\alpha_p}$ with $\alpha_i \in \mathbb {N}^*$ 
such that $gcd (\frac{\Pi -1} { d}, d) =1$.}
\medskip

By definition of $d$, we have $n_i = k_id +1$ where  $k_i$ are $p$ positive
coprime integers. So there exists $p$ integers $(\beta_1, ...., \beta_p)$ such
that  \ \ \ \ $k_1\beta_1+ ....+k_p \beta_p= 1$ \ \ \ (*).

\medskip

Fix $p$ positive integers  $\alpha'_i \in \mathbb {N}^*$  such that
$\alpha'_i>\vert \beta_i\vert $.

\medskip

If $k_1\alpha'_1+ ....+k_p \alpha'_p$ and  $d$ are coprime 
 then we set  $\alpha_i = \alpha'_i $.

\smallskip

If not,  $k_1\alpha'_1+ ....+k_p \alpha'_p = w'd'$  with $d'>0$ a divisor of
$d$. By multiplicating this equality by
$n=\frac{d} { d'}$,  we obtain :  $k_1 (n\alpha'_1)+ ....+k_p (n\alpha'_p) = w'd$.
 Adding (*) to this equality, we get  $k_1 (n\alpha'_1+\beta_1)+ ....+k_p
(n\alpha'_p+\beta_p) = w'd+1$ and  we set $ \alpha_i = n\alpha'_i+\beta_i >0$.
It's clear that $d$  is relatively  prime to  $w'd+1$.

\medskip

Hence in both cases,  we have found  $p$ positive integers $\alpha_i>0$ such
that  $w= k_1\alpha_1+ ....+k_p \alpha_p$  and $d$ are coprime. Consider  $\Pi =  n_1 ^{\alpha_1} n_2
^{\alpha_2}... n_p ^{\alpha_p}$  and compute $\Pi -1$. We denote by $C_n^p$ the
binomial coefficient $C_n^p= \frac{n!} { p! (n-p)!}$.
$$\Pi -1 = (k_1d +1) ^{\alpha_1} .... (k_pd +1) ^{\alpha_p} = 
\sum _{\begin{array}{ll}  i_j  \in \{0,...,\alpha_j\}\\   j \in \{1,...p\} \end{array}}  (C_{\alpha_1}
  ^{i_1}...C_{\alpha_p} ^{i_p})  (k_1 ^{i_1}....k_p ^{i_p}) d ^{
    ({i_1}+... {i_p})} \ \ \ - 1 = $$
$$
\sum _{\begin{array}{ll}    i_j  \in \{0,...,\alpha_j\} ,\sum i_j\geq 2  \\ j \in
  \{1,...p\}  \end{array} }  (C_{\alpha_1}
  ^{i_1}...C_{\alpha_p} ^{i_p})  (k_1 ^{i_1}....k_p ^{i_p}) d ^{
    ({i_1}+... {i_p})}  \ \ \ + \sum _{ \begin{array}{ll}  i_j  \in
      \{0,...,\alpha_j\}, \sum i_j \leq 1  \\
j \in   \{1,...p\}\end{array}} \hskip -1truecm  .....  \ \ \  \ \ \ - 1 = $$

$$= d^2 \Sigma   + d(k_1\alpha_1+ ....+k_p \alpha_p)= d^2 \Sigma   + dw{\text { \ where \ }} \Sigma \in \mathbb {N}.$$

Because,  $\sum i_j \leq 1 $ if and only if all $i_j$ are zero (in this case
the term of the sum is 1) , or  all $i_j$
are zero except one $i_{l}$ that is egal to $1$ (in this case
the term of the sum is $\alpha_l k_ld$). 

\medskip

 Using that  $w = \frac{\Pi -1 } { d} - d \Sigma $ and $d$ are  coprime,
we conclude that  $\frac{\Pi -1} { d} $ and $d$ are coprime.

\bigskip 

\noindent {\bf Step 2 : Construction of maps with irrational rotation numbers.}

\medskip

Fix $j \in \{1,....,p\}$ and set $r_k = k \frac{\Pi-1 } { d} $ for 
$k\in \mathbb{N}^*$, with $\Pi$ given by the previous step.

 Fix  $i \in \{1,....,p\}$ and  consider $f_i$ the
Boshernitzan on $S_{r_k}$ defined by its two slopes : 

\centerline {$\dis \lambda_1= n_i>1 {\text   { and }}  
 \lambda _2 =\frac{\lambda_1} { \Pi} <1.$}

\noindent For  proving that  $f$    is
in  $T_{r_k, (n_i)}$, it suffices  to prove that its break  $ a_i\in \mathbb
{Z}[\frac{1} { m}]$. We have $\dis a_i = r_k\frac{1-\lambda_1} {\lambda _2-\lambda _1 }=
k(\frac{\Pi-1 } { d}) \times \frac{1-\lambda _1} {\lambda _1 \Pi  ^{-1} -\lambda
  _1} = \frac{k \Pi } {\lambda_1 } \frac{\lambda_1 -1} { d}$. But,
$1-\lambda_1 \in dA$ so   $\frac{\lambda_1 -1} { d} \in A$, thus $a_i\in A$.

\medskip

$\underline{\text {To sum up}} $ : for any $i \in \{1,....,p\}$, we have construct 
$f_i \in T_{r_k, (n_i)}$ such that :  $$\rho (f_i) = \frac {\log n_i } {  \log \Pi }
\notin \mathbb{Q}  {\text { \ and \ }}  \sigma _{f_i} (0) = {\Pi}.$$

Futhermore, for $i \in \{2,....,p\}$, the $f_i$ have irrational rotation
numbers such that $\rho(f_i)$ and $1$ are  $\mathbb{Q}$-independent, 
  and commute because have the same jump $\sigma_{f_i}(0)= \Pi$.  Therefore  the Thompson-Stein group  $T_{r_k, (n_i)}$ contains a $\mathbb {Z}^
{p-1}$ that acts freely on $S^1$.

\medskip

The classes in $\frac{\mathbb{Z}} {d \mathbb Z}$ described by the  $r_k= k 
\frac{\Pi -1 }{d}$, when $k$ ranges over  $\mathbb{N}^*$ are the multiples of the class in
$\frac{\mathbb{Z} } { d \mathbb{Z}}$ of  $gcd( d, \frac{\Pi -1 } { d})=1$. 
So any 
positive integer $r$ is in the class $(mod \ d)$  of some $r_k$. By Bieri-Strebel
criterion,  the group  $T_{r, (n_i)}$   is PL-conjugate  to some  $T_{r_k, (n_i)}$,
thus it contains a $\mathbb {Z}^ {p-1}$ that acts freely on $S^1$. 
This ends the proof of theorem 2, part A.

\subsection{ Proof of theorem  2  B, C and C'. }

These  proofs  use a rigidity result for  PL-actions of $\mathbb {Z}^2$ on circles,
which   is a lemma in [18] :{\it   ``If  $f,g \in PL^+(S_r)$ generate  a
  freely acting group then $f$ and $g$ are PL 
conjugate to two  commuting Boshernitzan on $S_r$''}. Here, we slightly
modified the proof  given by Minakawa in order to obtain  further informations
on the conjugating homeomorphism  and to compute the rotation numbers  of $f$ and $g$.

\medskip

\noindent {\bf Definitions.}

- A homeomorphism $f\in PL^+(S_r)$ has the {\df (D)-property} if the
product of the $f$-jumps on each orbit is trivial. Or equivalently 
  the break set  $BP(f)$ is contained in the union of finitely
many pieces  of orbits $C_i = \{ a_i ,f(a_i ) ,  ..., f^l(a_i)\}$, $i \in I$ such
that  $\dis \prod_{c\in C_i} \sigma_f(c) =1$.  

- For $f$ with (D)-property,  we define
 $\dis \pi(f) =\prod_{i\in I, \ 0< k\leq l_i} \sigma_{f^{N+1}} (f^k(a_i)) $,
 where $N= max \  l_i $

 \vskip -0.3 truecm
\noindent   with $l_i$ defined as above.

\bigskip

\noindent {\bf Lemma 1.} {\it If  $f,g \in PL^+(S_r)$ generate a  freely acting
  group, then they  have the  (D)-property.}

\medskip

\noindent {\bf Proof of Lemma 1.}  Let $a$ be a break of $f$, we have to prove that the product of the $f$-jumps
on the $f$-orbit of $a$ is trivial. In fact it suffices to prove the
triviality of  the product of jumps on   finite pieces of orbit $\{f^{-M}(a),
....,f^{M}(a)\}$, for any  sufficiently large $M$.

\smallskip

Since $f$ has finitely many  break points  and since  the $f$-orbits
of the points $g^k (a)$ are disjoints, it's possible to find $k \in
\mathbb {N}$ such that the  $f$-orbit of $g^k(a)$ doesn't meet $BP(f)$. Thus, for all $n\in  \mathbb {Z}$, we have $\sigma_{f}(f^n (g^k(a))=1$. Futhermore, for all $x\in S_r$, we have :  $$\sigma_{g^k \circ f\circ g^{-k}
}(g^k(x))=  \frac {\sigma_{g^k} ( f(x)) \sigma_{ f} (x)} { \sigma_{g^k  }(x)}.$$

\noindent Using the commutativity hypothesis, we have 

 $$1= \sigma_{f}(f^n (g^k(a))= \sigma_{g^k \circ f\circ g^{-k}
}(f^n (g^k(a))= \sigma_{g^k \circ f\circ g^{-k}
}(g^k(f^n(a)) =   \frac{\sigma_{g^k} ( f(f^n(a)) \sigma_{ f} (f^n(a))} {
  \sigma_{g^k  }(f^n(a))}.$$

\noindent It follows that \ \ \ \ $\dis \sigma_{ f} (f^n(a)) = 
\frac {\sigma_{g^k} ( f^n(a))} { \sigma_{g^k  }(f^{n+1}(a))},$ \ \ \ hence the product

\medskip

\noindent $\displaystyle \prod_{n\in \mathbb {Z}}\sigma_f (f^n(a)) = 
 \prod_{n\in \{-M,..,M\}}\sigma_f (f^n(a)) $ \ \ \ \ telescopes to \ \ \ $
\frac{\sigma_{g^k} ( f^{-M}a))} { \sigma_{g^k  }(f^{M+1}(a))}$. 

\smallskip

\noindent Finally, by choosing $M$ sufficiently large (depending on $k$)
such that  $f^{-M}(a) $ and  $f^{M+1}(a)$ are  not breaks of $g^k$,  we get
that  the
product of the $f$-jumps along the $f$-orbit of $a$ is trivial. We conclude that $f$ (and also $g$)  has the (D)-property.

\bigskip

\noindent {\bf Lemma 2.}
 {\it Let   $f\in T_{r, \Lambda, A}$  with  (D)-property and irrational
   rotation number.  Then there exists a
 PL map $H_f$ that conjugates $f$ to the $S_r$ Boshernitzan example  ${\cal
   H}_r \circ  B^{\Pi(f)}_{\rho(f)} \circ   {\cal H}_r ^{-1}$. Furthermore :  

If $\Pi(f) \not=1$, then $\rho(f) =
\frac {\log \alpha } { \log   \Pi(f)}$, where $\alpha \in \Lambda$.  

If $\Pi(f) =1$,  $\rho(f) = \frac{ \alpha } { \beta}$, where $\alpha,  \beta
\in \Lambda$.

\noindent  The break set  $BP(H) \subset A$, the jump sets $\sigma (H )
\subset \Lambda$ and   $\sigma (B^{\Pi(f)}_{\rho(f)}) \subset \Lambda$. But,
in general,  the slope sets $\Lambda(H)  $ and  $\Lambda( B^{\Pi(f)}_{\rho(f)}) $ are  not contained in $\Lambda$.}

\medskip

\noindent {\bf Proof of Lemma 2.}  By definition of the (D)-property, the set
\ \ $\dis BP(f) \ \subset  \bigsqcup_{i\in I} {\cal
  C} _i$  \ with  

\vskip -0.3truecm 
\noindent  ${\cal C}_i  = \{ f^{k}(a_i),\  k=0, ...,l_i\}$  \ and  \ $\displaystyle
 \prod_{c\in  {\cal C}_i}\sigma_f(c) =1$.

\noindent  Consider a PL homeomorphism $H_f=H$ of $S_r$  such that  $f(a_i)$,
...., $f^{l_i}(a_i)$ for  $ i\in I$ are breaks of $H$ with associated
jumps $\sigma_H ({ f^{k}(a_i)) = \sigma_{ f^{N+1}} (f^{k}(a_i))}$
 for $k=1,...l_i$,  where   $N = max \ l_i$. Note that for $a_i$ we also have $\sigma_H ({ f^{k}(a_i)}) =
\sigma_{ f^{N+1}} (f^{k}(a_i))=1$.

\smallskip

\noindent A necessary and sufficient condition for $H$ to have  these breaks
and  no more is that the product of their $H$-jumps is trivial, that is  
$\Pi(f) =1$.

So, when $\Pi(f) =1$, we normalise $H$ by $H(0) =0$.

So, when $\Pi(f) \not=1$, we add a break point $c$ such that $c\in A
\setminus \cup \{a_i,..., f^{l_i}(a_i)\}$ and $\sigma_H (c) = (\Pi (f) )^{-1}$ and
normalise $H$ by $H(c)=0$.

\medskip

By definition, this map $H$ has breaks in $A$, jumps in $\Lambda$.
The slopes of $H$ are in $frac(A)=\{\frac{a} { b}  :  \ a,b\in A\}$ because  are 
ratio  of  lengths  of two intervals with endpoints in $BP(H)\subset A$.

\medskip

\noindent  \underline {Now we have :}

- $BP(H\circ f \circ H^{-1}) \subset \{H(a_i),..., H(f^{l_i}(a_i)), i\in I\} \cup \{
H(c), H(f^ {-1} (c) \}$,

\medskip

- for $i\in I$, $0\leq k \leq l_i$, the jump :

$$\sigma_{H\circ f \circ H^{-1}}(H(f^k(a_i))= \frac{\sigma_{H}(f^{k+1}(a_i)) \times 
\sigma_{f}(f^{k}(a_i))} {\sigma_{H}(f^{k}(a_i))}=
\frac {\sigma_{f^{N+1}}(f^{k+1}(a_i)) \times \sigma_{f}(f^{k}(a_i)) } { \sigma_{f^{N+1}}(f^{k}(a_i))}=1,$$

- the jump  $\sigma_{H\circ f \circ H^{-1}}(H(c)) = \Pi(f)$,

- the jump  $\sigma_{H\circ f \circ H^{-1}} (H(f^ {-1} (c) ) = \Pi(f)^{-1}$.

\bigskip

\noindent  Thus, when $\Pi(f) \not=1$, the map $F= H\circ f \circ H^{-1}$ has exactly two
breaks $0=H(c)$ and $F^ {-1} (0) =H(f^ {-1} (c))$, it's a  Boshernitzan on
$S_r$ ; its jump at the initial break $0$ is $\Pi(f)$ and 
its rotation number is $\rho(f)$ but also \ 
$\dis \frac{\log DF_+(0)} { \log DF_+(0)-\log DF_-(0)}= \frac{\log DF_+(0)} {
  \log\Pi(f)}$ \ as it's a Boshernitzan. 

\smallskip

\noindent Moreover,    $\dis  DF_+(0)= D_+(H\circ f \circ H^{-1}) (H(c)) =
\frac {DH
  (f(c))  Df (c)} {  DH_+(c)} =  Df (c)\frac {DH (f(c))} { DH_+(c) } \in \Lambda$, as $Df (c)\in
\Lambda$ and the ratio of the derivatives of $H$  at two points 
 is  a product of $H$-jumps.

\medskip 

So $\dis \rho(f) =\rho(F) = \frac{\log \alpha } {\log \Pi (f)}$, with $\alpha$
and  $\Pi(f)$ in $\Lambda$.

\bigskip

\noindent  Now, when $\Pi(f) =1$, the map $F= H\circ f \circ H^{-1}$ 
is PL and has no 
break so it's a rotation of angle $\rho(f)$ but also of angle  $F(0)= 
H\circ f \circ H^{-1} (0) = H(f(0))= DH({f(0)}).  f(0) + \beta \in frac A$, as  $DH(f(0)
\in frac A $ and $\beta \in frac A$ (this can easily be etablished 
by writing that $H$ is continuous at its breaks $\alpha_i \in A$, that $H(0)=0$
and that the slopes are in $frac A$).  This ends the proof of Lemma 2.

\medskip

\noindent {\bf Remarks on Lemma 2.}

- When $\Lambda$ is a  subset of
  $\mathbb{Q}$, the case $\Pi(f)=1$ can't occur, 
since $\rho(f) $ would be rational.

- If $f$ and $g$  satisfy the hypothesis of lemma 2 and commute then $\pi(f)
=\pi(g)$, since the conjugation to rotation with $h(0)=0$ is $R_{H^{-1}_f(0)}\circ H_f
\circ h_{\pi(f)} =R_{H^{-1}_g(0)}\circ H_g \circ h_{\pi(g)}$ and an easy calculation gives
that $ h_{\pi_1}  \circ h_{\pi_2}^{-1}$ is PL if and only if $\pi_1=\pi_2$.

- For similar reasons, if  $f$ and $g$  satisfy the hypothesis of lemma 2 and
are PL conjugate  then $\pi(f)=\pi(g)$. 

\bigskip

\noindent {\bf Proof of Theorem 2.B.}   By absurb, we suppose that there is
$\mathbb {Z}^p\subset  T_{r,  (n_1,...,n_p)}$   acting  freely.

If $p=1$, this is impossible by   theorem 1 A.

If $ p\geq 2$,  there exist $p$ commuting homeomorphisms $f_i$ of 
irrational rotation numbers  $\rho(f_i)$  affinely $\mathbb {Q}$-independent. Using  lemma 1 and
2 (and its remarks) we get  $\rho(f_i) =\frac {\log \alpha_i } { \log \Pi(f_i)}$.

Furthermore all $\Pi(f_i)$ are egal to a same $\Pi$ (the $f_i$ commute) and
  the  $\rho(f_i)= \frac{\log \alpha_i } { \log \Pi}$ and $1$ must be 
 $\mathbb{Q}$-independent. This means (multiplicating by $\log \Pi \in\mathbb{Q}$)
 that the $ {\log \alpha_i }$ and $\log \Pi$ are $\mathbb {Q}$-independent. 
But, $ {\log \alpha_i }$ and $\log \Pi$  are $(p+1)$ real numbers 
in the additve group generated by the $p$ real numbers $ \log n_i$, so
they are $\mathbb{Q}$-dependent, a contradiction.

\bigskip

\noindent {\bf  Proof of Theorem 2.C.} If  $f\in T_{r, (n_i)}$ is contained in a $\mathbb{Z}^2$ acting freely, by
lemma 1, the map $f$ has  (D)-property. Thus, by lemma 2 (and its remarks),
the  rotation number of $f$ is  $\frac{\log \alpha } { \log  \Pi(f)}$
 with $\alpha $ and $ \Pi(f)$ in $<n_i>$. This proves the claim.

\bigskip

\noindent {\bf  Proof of Theorem 2.C'.} Fix $p\geq 1$ and  fix
 $\rho =\frac {\log \alpha } {\log \beta} \in  ]0,1[$  with $\alpha$, $ \beta
$ in $\Lambda =<n_i>$ and $1<\alpha < \beta $. As ${ \beta}-1 \in (1-
\Lambda)A =(d\mathbb{Z}) \Lambda$ we can write  \ \ \   ${ \beta}-1 = d  n_\beta\lambda$ with $ n_\beta\in \mathbb{Z}$ and
$ \lambda\in
  \Lambda$.

\medskip

Let $ r = d n_{\beta}$ and consider  the
Boshernitzan $f$ on $S_{r}$ defined by its two slopes 

\centerline {$ \dis \lambda_1=\alpha  {\text
  { and }}   \lambda _2 =\frac{\lambda_1} { \beta}= \alpha \beta ^{-1} <1.$ }

For  proving that  $f$    is
in  $T_{r, (n_i)}$, it suffices to prove that its break  $ a_i\in \mathbb
{Z}[\frac{1} { m}]$. We have $\dis a_i = r\frac {1-\lambda_1} {\lambda
  _2-\lambda _1 }= d n_{\beta}  \times \frac 
{\alpha -1} {\alpha \beta ^{-1}  ( {\beta} -1) }
= \frac{\alpha -1 } { \lambda\alpha \beta ^{-1} } \in A $. The  rotation number
of $f$ is $\rho(f) = \frac{\log \lambda_1} {  \log\lambda_2
-  \log\lambda_1}=
{\frac{\log \alpha } { \log \beta}}$. But $r=dn_\beta$ equal  $d \ (mod \ d)$ 
so using Bieri-Strebel criterion the group $T_{d, (n_i)}$ contains  a
PL-conjugate to $f$ that is a map of rotation number  $\frac{\log \alpha } {\log \beta}$.

\medskip

Remark also that  for  $p=1$, the result is  the consequence 2 of theorem 1.
\subsection {Proof of the consequences of theorem 2.}

Consequence 1 is a reformulation  of  A and B. The consequences 2, 3 and
4 use easy technical lemmas (lemma 4 and 5) and  a {\bf  Brin version ([4])
  of  the Rubin theorem of [21]  :} 

{\it Let $G$ be a
  group and  $\phi$,    $ \varphi$ be two representations  of $G$ in
  Homeo${}^+ (S^1)$. If the subgroups  $\phi (G)$ and  $ \varphi(G)$ are both locally dense then 
 $\phi$ and  $ \varphi$ are  topologically conjugate (that is there exists $h \in $ 
 Homeo${} (S^1)$ such that $\phi (g) = h  \circ  \varphi(g)\circ h ^{-1}$,
 $\forall g \in G$).

\smallskip

\noindent Where a subgroup $H$ of  Homeo${}^+ (S^1)$  is said to be {\df locally dense}  if for all 
open subset $U$ of $S^1$, the group $G_U = \{ g \in H : supp \ g\subset U\}$ has
every orbit in $U$ locally dense (i.e dense in some non empty  open subset of $U$).}

\medskip

\noindent {\bf Lemma 3.} {\it Fix two  integers $r\geq 1$ and $ k\geq 2$, then  the group
  $T_{r,k}$ is a  locally dense subgroup  of  Homeo${}^+ (S^1)$.}

\medskip

\noindent {\bf Proof of  lemma 3.} Let $U$ be an open subset of $S^1$
 and  $x_0\in U$. We have to prove that the closure $cl(G_U (x_0))$ of
 $G_{U}(x_0)$ contains
 an open subset of $U$. Let's choose $a_0, b_0\in A= \mathbb{Z}[\frac{1} { k}]$ such that $x_0 \in I_0= ]a_0,b_0[\subset U$. We'll prove that   $cl(G_{I_0} (x_0))$  contains a non empty open subset.

\medskip

\noindent For $\alpha\in ]0,(k-1)min ({x_0-a_0},\frac  { b_0 -x_0} { k})[ \
\cap A$, we consider the PL-homeomorphism $f_\alpha$ defined by  :

\medskip

$\dis f_\alpha(x) = kx + (1-k)a_0 $ \ \ \  \ if  \ \ $\dis x \in [a_0,
a'_0]$  \ with  $a'_0= a_0 +\frac {\alpha } { k-1}$,

$\dis f_\alpha(x) = x + \alpha  $   \ \ \ \  \ \ \ \ \ \  \ \ \  \ \ \ if \ \
$\dis x \in [ a'_0 , b'_0]$ \ with $b'_0=  b_0 - \frac{k\alpha } { k-1}$,

$\dis f_\alpha(x) = \frac{1} { k}x + ( \frac {k-1} { k }) b_0 $ \ \ \ 
if  \ \ $\dis x \in [ b'_0, b_0 ]$,

$\dis f_\alpha(x) = x$ \ \ \ elsewhere.

\vskip  - 5 mm 

\hskip 5,3 truecm \begin{picture}(100,100)

\put(0,0){\line(1,0){90}}  
\put(0,90){\line(1,0){90}}    
\put(0,0){\line(0,1){90}}    
\put(90,0){\line(0,1){90}}  

\put(30,-1){$\vert$}
\put(28,-10){$a_0$}

\put(80,-1){$\vert$}
\put(80,-10){$b_0$}

\put(-2,30){$-$}
\put(-10,30){$a_0 $}

\put(-2,80){$-$}
\put(-10,80){$b_0$}

\put(30,30){${\scriptscriptstyle \bullet}$}
\put(0,0){\line(1,1){30}}
\put(30,30){\line(1,2){10}}

\put(40,50){\line(1,1){20}}
\put(40,50){${\scriptscriptstyle\bullet}$}

\put(60,70){${\scriptscriptstyle\bullet}$}
\put(60,70){\line(2,1){20}}

\put(80,80){${\scriptscriptstyle \bullet}$}
\put(80,80){\line(1,1){10}}

\put(40,-1){$\vert$}
\put(40,-10){$a'_0$}

\put(60,-1){$\vert$}
\put(60,-10){$b'_0$}

\put(15,8){${1}$}
\put(50,63){${1}$}
\put(38,38){${k}$}
\put(68,82){$\frac {1} { k}$}
   \end{picture}

\vskip 0.4 truecm

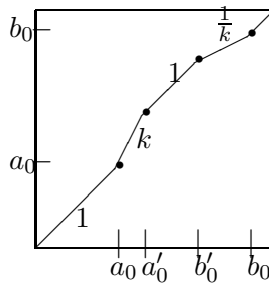
\centerline {\it Fig. 2}

\smallskip

We have $f_\alpha \in T_{r,k}$, its support  $supp f_\alpha = [a'_0,b'_0] \subset
I_0\subset U$  and  $f_\alpha (x_0) = x_0+\alpha$    for all 
$\alpha $ in the dense set $ A'= ]0,(k-1)  min ({x_0-a_0}, \frac{ b_0 -x_0} 
  {k})[ \ \  \cap  \ A \ $ as $a_0\leq
a'_0 <x_0 <b'_0 \leq b _0$. 

  Thus, cl$(\{f_\alpha(x_0), \alpha\in A'\})$
has non empty interior and is contained in the interior of $cl(G_{I_0}
(x_0))$. This ends the proof of lemma 3.

\bigskip

\noindent {\df \underline {Consequences} :  The groups $ T_{r, (n_i)}$ are locally dense
  subgroups} (since they contain some $ T_{r, k}$).

 {\df \hskip 2 truecm Each isomorphism between such groups  and each 
   automorphism  of such groups  is 

\hskip 2 truecm realized by a conjugation  in  Homeo $(S^1)$.}

\bigskip

{\bf Lemma 4 : } {\it   Let $\Gamma_1=<f_1, g_1> \  \subset $  Homeo $ ^+(S_{r_1})$ and
  $\Gamma_2=<f_2, g_2> \ \subset $ Homeo $ ^+(S_{r_2})$  be two $\mathbb {Z}^2$  that 
act freely on circles and such that $\rho (f_1) = \rho (f_2)$ and   $\rho (g_1) = \rho
(g_2)$. Then  $\Gamma_1$ and $\Gamma_2$ are conjugate, and the maps such  that  
 $\Gamma_2 = h\circ \Gamma_1\circ h^{-1}$ are  

\smallskip
$\dis h_b = H_{f_2} ^{-1} {\circ {\cal H} _{r_2}} \circ  h_{\pi(f_2)} \circ R_b
\circ   h_{\pi(f_1)} ^{-1} {\circ {\cal H} _{r_1}} ^{-1}\circ H_{f_1} \ \ 
{ \text  { with  } } b\in S^1, $

\smallskip

\noindent where $H_f$ is the map defined in lemma 2, $ {\cal H} _{r}$ is the homothety of
ratio $r$ and  $h_{\pi}$ is the conjugation between the Boshernitzan $B_\rho
^\pi$ and the rotation by $\rho$. }

\medskip

\noindent The proof of  lemma 4 is the combination of lemma 2 and properties of 
 Boshernitzan-examples and we leave it to the reader.

\medskip

{\bf Lemma 5 : }{\it  Let $\pi_1$, $\pi_2$ in $]0,1]$ and $a,b$ in $S^1$,
  denote   $F_b = h_{\pi_2}\circ R_b \circ h_{\pi_1}^{-1}$, \break then
$F_b\circ R_a \circ F_b^{-1} $ is a PL-map if and only if $a=0$ or $\pi_1=1$ or
 $\pi_1= \pi_2$. }

\medskip

{\bf Proof of lemma  5.}  The  calculus are done in  $[0,1[$, circle
homeomorphisms $f$ being identified with  bijections $\tilde f$ (mod 1)  of  $[0,1[$.

\medskip

We write  $h_\pi (x) =\frac { \pi ^x -1 } { \pi -1}=  \frac{ 1 } { \pi -1}(e^{x\log
  \pi}  -1)$ et $h_\pi^{-1} (x) =\frac{1} 
 { \log \pi} \log ( x(\pi -1) +1)$.
\medskip

We compute  \ \ \ $\displaystyle F_b(x)=  h_{\pi_1} \circ R_b\circ   h_{\pi_2
}^{-1} (x)= \frac{1}{  \pi_1 -1} \left ( \ \pi_1^{b(x)}(x( \pi_2 -1) +1 )
  ^{\frac {p} {s}}-1 \
 \right),$ where  $p= \log \pi_1$, \ $s =\log \pi_2$  \ and \ $b(x)\in \{b, b-1\}$.

\medskip
Finally, $\dis F_b\circ R_{a}\circ F_b^{-1}(x) =\frac {1} {  \pi_2 -1}  \left ( \pi_2^{-b_1(x)}
\left ( \  a (\pi_1-1) + 
\pi_1 ^{b_2(x)}(x( \pi_2 -1) +1 )^{\frac {p}{ s  }} \right) ^{\frac {s}{ p} }-1 \
\right)$,  where $b_1(x),b_2(x)\in \{b, b-1\}$.
 
\medskip

Finally, $ F_b\circ R_{a}\circ F_b^{-1}(x) $ is PL   if and only if
$a'(\pi_1-1)= 0$  or  $\frac {p} { s} =1$, this gives lemma 5.

\bigskip

\noindent  {\large \bf Proof of  consequence 2 : } `` Topological rigidity implies PL rigidity''.

\medskip

Let  $\phi$ be a  representation  of $G= T_{r, (n_i)}$, $p\geq 3$  in   $PL^+
(S_{1})$  that  is topologically conjugate by $h$ to the standard one. By
eventually composing $h$ with the standard involution
($x\mapsto 1-x$), we may assume that $h$ is orientation preserving. 
 Consider a rank 2 subgroup $<f_1,g_1>$ contained in the   ``free'' $\mathbb {Z}
^{p-1}$ subgroup of $G$ that has been   constructed in the  previous section
(proof of theorem 2.A). Changing $r$ in
 its class mod $d$, we may assume that $f_1$ is a Boshernizan : $f_1 = 
{\cal H}_{r_1} \circ h_{\pi_1}\circ R_{\rho_1} \circ h_{\pi_1}\circ {\cal
  H}_{r_1}$.
\medskip

 The image of   $<f_1,g_1>$  is a  $\mathbb{Z} ^{2}$ acting freely,  we can apply
 lemma 4 for  $<f_1,g_1> \subset  G$ and $<\phi(f_1),\phi(g_1)> \subset \phi(
G)$. We obtain that  $h$ is one of the $h_b$. Let $R_a\in G$ be a non trivial rotation, we have  $\phi(R_a) =
h_b \circ R_a \circ h_b ^{-1}\in PL^+ (S^1)$. Replacing $h_b$ by its
expression, it means 

\smallskip

\noindent that   $\dis H_{f_2} ^{-1}  \circ  h_{\pi(f_2)} \circ R_b
\circ   h_{\pi(f_1)} ^{-1} \circ {\cal H} _{r_1} ^{-1}\circ R_a \circ 
 {\cal H} _{r_1} \circ h_{\pi(f_1)}\circ R_b^{-1}  \circ  h_{\pi(f_2)}^{-1}
 \circ \  H_{f_2} \in PL^+ (S^1).$

\smallskip

Using the fact that $ H_{f}$ are PL-maps and lemma 5,  we get that ${\pi(f_1)} 
=1$ that  is impossible as   $f_i$ is a rational homeomorphism with irrational 
rotation number or  ${\pi(f_1)}={\pi(f_2)} $ in this case $ h_{\pi(f_2)}
\circ R_b \circ   h_{\pi(f_1)} ^{-1} $ is a Boshernitzan so it's PL. Finally, we
have proved that $h_b$ is PL.

\medskip

\noindent  {\large \bf Proof of  consequence 3 : } `` Non exoticity of automorphisms and isomorphisms''.

According to  Rubin theorem,  isomorphisms and  automorphisms of
Thompson-Stein groups are realized by topological conjugations.
 According to  the previous point,  they are realized by   PL-conjugation,
 provided the fact that $p\geq 3$. More precisely, it's possible to give a
 complete  

\smallskip

\noindent {\df  description  of  PL-automorphisms.} This 
 can be found in [1] or in [16] and prove that a PL
automorphism of $G=T_{r, \Lambda, A}$ has slopes in $inv A$ and breaks in
$A$.

\medskip

\noindent {\large \bf Proof of consequence 4 :}  ``Distincts  Thompson-Stein groups are
  not  isomorphic''.

Suppose, by absurd, that  $G=  T_{r, (n_i)}$ and $G'=T_{r, (n'_j)}$ are
isomorphic and $G$ is of free rank at least 2.  By consequence 3, the groups  $G$ and $G'$
are PL-conjugate. Using the summary at the end of proof of theorem 2A, for all
$i\in \{1,..,p\}$ the group $G=T_{r, (n_i)}$ contains a map $f_i$ with
rotation number $\rho(f_i) =\frac {\log n_i } { \log \Pi}$. The image  $f'_i$ of
$f_i$ in  $G'=T_{r, (n'_i)}$  is contained in a free $\mathbb {Z}^2$ 
 so   its  rotation number  $\rho(f'_i)=\frac {\log \alpha' } {  \log \Pi(f'_i)}$,
$\alpha'\in <n'_j>$ and as  $f'_i$ is PL conjugate to $f_i$ we also have 
$ \Pi(f'_i) =  \Pi$  (according to  remarks on lemma 2). So, for all $i\in \{1,...,p\}$ the integer $n_i\in <n'_j>$ and conversely.
 Consequently, the groups  $ <n_i> =<n'_j>$ and  eventually changing the basis
 of the slope sets we  have $n_i= n'_j$,  hence the groups $G=G'$.

\section {Proof of theorem 3.} 
\subsection {Dynamical properties for representations of  $\mathbf{G} =
  \mathbf { T_ {1, ( n_i  ) }}.$} In this section, we study  the dynamics 
of a non trivial  representation   $\dis \phi : G \rightarrow  Homeo ^+( S^1)$.  

\medskip

{\bf Notations.} $ F= \{g\in G : g(0)=0\}$,  

\hskip 2 truecm  $ F'=[F,F]= \{g\in G : g(0)=0  \ {\text {and }}\  Dg(0) =1\}$,

\hskip 2 truecm  $A = \mathbb {Z} [\frac {1} { m}]$, where  $m = lcm (n_i)$ is identified to  the  subgroup of $G$ consisting
 of  rotations  i.e  the  group of rotations of angles in $A$.

 \bigskip

\noindent {\bf Algebraic properties, useful for  dynamical properties on  $F$
  and  $G$ :}

\smallskip

-  $F$ doesn't  contain any free group of rank 2 (Brin-Squier) and has no 
 torsion.

-  $F'$  is  simple  (Bieri-Strebel, Stein).

-  $T_{1,k}$, for all  $k\geq 2$ integer in $<n_i>$, 
is  simple (Brown),  is contained in $G$ and contains $A$.

\bigskip
\noindent {\large \bf Proposition 1.} {\bf  The  group $\phi(F)$ acts on
  $S^1$ with a  finite orbit.}

\medskip

 { \bf Proof of proposition 1.} By  a Margulis theorem [Ma 2000],  since
 according to Brin-Squier  $F$
 contains no free group of
rank 2, the group $\phi(F)$  preserves  a  measure $\mu$ on  $S^1$. Thus,  the 
 rotation number map  $\rho$ is a representation  of  $F$
into an abelian group. In particular, it's trivial on the first derivated
group $\phi([F,F])$ and consequently, any element of $\phi([F,F])$ has a fixed point.

\bigskip

\medskip

{\it The  {\df dynamical alternatives theorem} (see [Gh2000] prop. 5.6)  says that
$\phi(F)$ satisfies one of the mutually exclusive following possibilities  :

\hskip 1 truecm 1. there  is  a finite orbit  or

\noindent there  is  a  unique  minimal $K$ that is   the accumulation set
of any orbit and 

\hskip 1 truecm 2. $K=S^1$ and all orbits are   dense or 

\hskip 1 truecm 3. $K\not=S^1$ is said to be exceptionnal ($K$ is homeomorphic
to the  Cantor set).}

\medskip

If  $\phi(F)$  has no    finite orbit, the support $supp \mu$ of $\mu$
coincide with  the unique  minimal  $K$ of $\phi([F,F])$. If $K= S^1$ then  by
setting  $h(x) = \mu([0,x])$,  one constructs a conjugation from  $\phi(F)$ to
a  group of   rotations and $F$ is an abelian group, this is not the case. 

So, $K=supp \mu $ is a cantor set and  any point in $K$ is fixed by
 $\phi([F,F])$. 

Moreover,  let  $g_i$, $i=1,..;p$  be the  generators of  $F$, they have
rational  rotation  numbers ${p_i\over q_i} $, since if not by  
 Denjoy et unique ergodicity  $\mu$ would be of total  support. 

Then the  finite type abelian group    $H ={ F \over [F,F]}$ acts on $K$
(since $[F,F]\subset Ker \phi_{\vert K}$) and any  generator $g_i$  is of 
finite order $q_i$. So  for all  $x\in K$, we have  $H(x) =
 \{g_1^{s_1}.....g_p^{s_p} (x), 0\leq  s_i <q_i\} $  that is finite.

Finally, the $\phi(F)$ orbit of any $x\in K$  is   the set of the 
$\phi(f). x_0 = \gamma c (x_0)$, where $ c\in  \phi([F,F])$, $\gamma \in
\{g_1^{s_1}.....g_p^{s_p},   0\leq  s_i <q_i\} $, so it's finite since $
\gamma c (x_0) = \gamma(x_0)$. This contadicts the dynamical alternatives 
theorem.  This ends  the proof of the  proposition 1.

\bigskip

\medskip

{\bf Proposition 2.} {\it The groups 
$\phi (G)$, $\phi(A)$ and  $\phi(T_{1,k})$ for  $k $ integer in $ <n_i>$  have no  finite
orbits   and have a unique same minimal $ K$ that is 
  the accumulation set of  any orbit.}

 \bigskip

{\bf Proof of  Proposition 2.}

\medskip

\noindent {\bf Step 1 : the groups $\phi (G)$, $\phi(A)$, $\phi(T_{1,k})$ for
  $k$  integer in $<n_i>$ have no finite orbits.}
 
\smallskip

The group  $\dis \phi( \mathbb {Z}[{1/k}])$  is isomorphic to $ \dis\mathbb {Z}[{1/
  k}]$    by proposition 2,  so it contains element
of arbitrary big order, so it has no finite orbits. Now, the groups $\phi(G)$
and   $\phi(T_{1,k})$  contain some $ \dis \phi(\mathbb {Z}[{1/ k}])$, so they  also
have no  finite orbits.

\medskip

\noindent {\bf Step 2 : the groups $\phi (G)$, $\phi(A)$, $\phi(T_{1,k})$, for
  $k$ integer in   $<n_i>$  have the   same unique  minimal $ K$.}

\smallskip

\noindent According to dynamical alternatives theorem, each of  
these groups $\phi (H)$ have a unique minimal denote by $K(H)$.  
 From inclusions between these groups, follow the  inclusions :

\smallskip

$\dis K(  \mathbb {Z}[{1 /    k}]) \subset K(T_{1,k}) \subset K(G),$  $k$ integer in $<n_i>$ 

$ \dis K(  \mathbb {Z}[{1 /    k}] )\subset K(A) \subset K(G) .$

\medskip

\noindent{\bf  We first prove that $K(A)= K( \mathbb {Z}[{1 /    k}])$.}  \ \ --For
any $a\in A$ and any   $\dis  \kappa \in \mathbb {Z}[{1/ k}]$,   we have \ \ \ 
 $\dis R_\kappa\circ R_a ( K( \mathbb {Z}[{1/ k}]))= 
R_a\circ R_\kappa ( K( \mathbb {Z}[{1/ k}])) =R_a( K( \mathbb {Z}[{1/k}]))$. 
This means that the set $\dis R_a ( K( \mathbb {Z}[{1/ k}]))$ is  $\dis \mathbb
{Z}[{1 /    k}]$-invariant, so by minimality and unicity of $K( \mathbb {Z}[{1/k}])$ this
set contains $\dis K(\mathbb {Z}[{1 /    k}])$. 

\smallskip

As the  same holds for $(-a)\in A$, we get  $\dis K( \mathbb {Z}[{1/k     }])$ is
$A$-invariant, so     $\dis K( \mathbb {Z}[{1/     k}])$  contains $K(A)$ by
minimality and unicity, this is the missing inclusion.-- 

\bigskip

\noindent{\bf  Now, we prove that $K(A) = K(G)$.} \ --As each  $g\in G$ can be written as $ g= R_a \circ f$, 
with $f\in F$, $a \in A$ ;  we have  $ \phi(g) = \phi(a) \circ \phi(f)$. 
But we have proved that there exists a 
point $x_0$ with a finite  $\phi (F)$-orbit : $\{x_0, ..., x_p\}$, so  
$ \phi(g).x_0 \in  \bigcup_{i=0}^p   \phi(a) . x_i$, for all $g \in G$.
 Thus  $ \phi (G) .x_0 =  \bigcup_{i=0}^p   \phi(A) .x_0$ and so on for the 
accumulation sets that are exactly the $ K(\ )$.

 Hence, we have $K(G) = K(A)$.--
\medskip

Finally, we have proved that  $K(\mathbb {Z}[{1 /    k}]) = K(A) = K(G)$ and
using the first inclusions we have that this set is also $
K(T_{1,k})$. According to the dynamical alternatives theorem, this set is
the accumulation set  of any orbit.

\subsection{\bf Proof of theorem 3.A, from now $\mathbf {G =  T_{
    1, ( 2,   { n}_2,..., n_p )}, p\geq 2}$.}

{\bf Step 1.} Any non trivial representation $\phi$ of $G$ in PL${}^+
(S^1)$ or Diff$^2_+(S^1)$ such that $\phi(G)$ has  
dense orbits is conjugate to the standard PL  representation.

\medskip

{\df Proof} -- The induced action on the 
classical Thompson group $T_{1,2}\subset G$ satisfies the hypothesis of
III proposition 3.9  p 230 in [11]  : it's non trivial  and the group $\phi
(\mathbb {Z}[{1/2}]) $  has dense orbits
because it has the same minimal than $\phi(G)$. Thus,  by Ghys-Sergiescu
results,  $\phi(T_{1,2})$ is conjugate to the PL standard
action $T_{1,2}$. So both actions $G$ and  $\phi(G)$ satisfy the hypothesis of
Rubin-Brin result, they are topologically conjugate.--

\bigskip

{\bf Step 2.}  Any non trivial representation $\phi$ of $G$ in Homeo${}^+
(S^1)$ is semi-conjugate to a representation  
$\phi_1$ of $G$ in Homeo${}^+(S^1)$ such that $\phi_1 (G)$ has  dense orbits.

\medskip

{\df Proof} -- If $K(G)$ is $S^1$, the representation $\phi(G)$ has  dense
orbits.  If $K(G)$ is not $S^1$, it's a Cantor set, the devil stairs associated
 to $K$ semi-conjugates $\phi$ to $\phi_1$ with  dense orbits.--

\medskip

\medskip

It's noted in [11]--proof of III Theorem 3.12 p 233-- that  when  $\phi$ is $C^2$
or $PL$,  the action  $\phi_1$ failed not to be
so regular but it satisfies all thinks that are needed for proving that
$\phi_1(T_{1,2})$ is conjugate to the PL standard action $T_{1,2}$ which implies
in our case that both actions $G$ and  $\phi_1(G)$ satisfy the hypothesis of
Rubin-Brin result, and hence  are topologically conjugate. We are now able to give 
 
\medskip

\noindent {\df the end  of  the proof of theorem 3.A.} Let $\phi$ be a
representation of $G$ in $PL(S^1)$ or Diff $(S^1)$, composing $\phi$ with the
standard involution ($x\mapsto 1-x$) we can suppose that $\phi $ is
orientation  preserving.  If   $\phi(G)$ has  
dense orbits then it is conjugate to the standard PL  representation, by step
1. If not, by step 2 and the previous remark,  $\phi$ is semi-conjugate to
such an action $\phi_1$, so it's semi-conjugate  to the standard PL action. But,
as $p\geq 2$, the group $G$ contains at least one element of irrational rotation number for which Denjoy theorem applies, this
and   unique ergodicity  implies that   the semi-conjugation is a conjugation.

\subsection{\bf Proof of theorem 3.B.}

Consider the Thompson-Stein group $G= T_{1, (2,n_2,..n_p)}$ with $p\geq
2$. By the theorem 2C', the group $G$ contains an
homeomorphism $B$ of  irrational rotation number  $\alpha =\frac{\log 2} { \log n_2}$  a diophantine number (see [21]) of  irrationality measure $\beta$   (for example for $n_2=3$ we have $\beta =7,616...$.)

Suppose that   there exists a representation $\phi : G
\rightarrow$ Diff $ {} ^k(S^1)$, for some $k\geq 2$. By theorem 3.A, 
$\phi (G)$ is  topologically conjugate by some homeomorphism $h$ 
to  $G$.

\smallskip

So $\phi  (B)= h \circ B \circ h^{-1}$ is a   $C^k$-diffeomorphism  with 
diophantine rotation number $\alpha$  of irrationality  measure $\beta$, by
Katznelson-Ornstein  result ([15]) it's $C^1$-conjugate to 
$R_\alpha$ provided $ k \geq \beta+1$.

\smallskip

 Thus for  $ k \geq \beta+1$, there
exists a $C^1$-diffeomorphism $H$ such that  $ R_\alpha = H \circ
h\circ B\circ h^{-1} \circ H^{-1}$. As $B$ is $PL$-conjugate to a Boshernitzan of rotation number $\alpha$, it is 
piecewise $C^\infty$ conjugate to $R_\alpha$. Using the unique ergodicty of
$R_\alpha$, we get that the map $H \circ h$ is piecewise
$C^\infty$, so  we get that   $h$ is piecewise $C^1$ with a finite number $N$ of breaks.

\smallskip

Finally, every  $g\in G$ can be written as 
 $ g = h ^{-1} \circ \phi(g) \circ h$ so it is piecewise
$C^1$ and has  at most  $2 BP (h) = 2N$ breaks, this is impossible as $G$
contains elements with arbitrary large  number of breaks.

 \subsection{\bf Proof of the consequence of theorem 3.} This is an immediate
 consequence of theorem 3A and  consequence 2 of theorem 2 (``topological
 rigidity implies PL rigidity'').

\eject

  \section{ References.}

\medskip

\smallskip


\hskip -1 truecm 1. \     R. Bieri and R. Strebel, {\it On groups of PL
  homeomorphisms of the real line.} Preprint, Math. Sem. der Univ. Frankfurt, Frankfurt am Main (1985).

\smallskip

\hskip -1 truecm 2. \  M. Boshernitzan, {\it Dense orbits of rationals.}  
 Proc. AMS, {\bf 117} (4), 1201-1203 (1993).

\smallskip

\hskip -1 truecm 3. \  M. G. Brin, {\it The chameleon groups of Richard J. Thompson:
 Automorphisms and dynamics},  Inst. Hautes \'Etudes Sci. Publ. Math. {\bf
 84}, 5--33  (1996).

\smallskip


\hskip -1 truecm 4. \ M. G. Brin, {\it Higher dimensionnal Thompson's groups.}
 Geom. Dedicata  {\bf 108}, 163--192 (2004).

\smallskip


\hskip -1 truecm 5. \ M. G. Brin, C. Squier, {\it  Groups of piecewise linear
  homeomorphisms of the real line.}  Invent. Math.  {\bf 79} (3), 
 485--498 (1985).

\smallskip


 \hskip -1 truecm  6. \  M. G. Brin,  F.  Guzmán, {\it Automorphisms of
 generalized Thompson groups.} J. Algebra {\bf 203}(1), 285--348 (1998).

\smallskip


\hskip -1 truecm 7. \ K. S. Brown, {\it Finiteness properties of groups.}  J. Pure
 Appl. Algebra {\bf  44}, 45--75. (1987).
\smallskip


\hskip -1 truecm 8. \  J. Cannon, W. Floyd, W. Parry,  {\it Introductory notes on
 Richard Thompson's groups.}  L'Ens. Math.  {\bf 42},  215-256 (1996).

\smallskip


\hskip -1 truecm 9. \ A. Denjoy, {\it Sur les courbes d\'efinies par  les \'equations diff\'eren\-tielles
\`a la surface du tore.}  J. Math. Pures Appl. {\bf 11}, 333-375 (1932).

\smallskip


\hskip -1 truecm 10. \  D.S Farley, {\it Proper isometric actions of Thompson's groups on Hilbert
  space.}  Int. Math. Res. Not. {\bf 45},  2409--2414 (2003). 

\smallskip


\hskip -1 truecm 11. \  \'E. Ghys, V. Sergiescu, {\it Sur un groupe remarquable de
  diff\'eo\-morphismes  du cercle.}  Comm. Math. Helv. {\bf 62}, 185-239
 (1987).

\smallskip


\hskip -1 truecm 12. \  \'E. Ghys, {\it Groups acting on the circle.}
L'Ens. Math.  {\bf 62}, 215-256 (2000).

\smallskip


\hskip -1 truecm 13.  \  M. Herman, {\it Sur la conjugaison diff\'erentiable des diff\'eo\-morphismes du cercle
\`a des rotations.} Inst. Hautes Etudes Sci. Publ. Math. {\bf 49}, 5-234 (1979).

\smallskip


\hskip -1 truecm  14. \  G. Higman, {\it  Finitely presented infinite simple
  groups.} Notes on Pure Mathematics {\bf  8}  (1974). Department of Pure
Mathematics,Department of Mathematics, I.A.S. Australian National University,
Canberra, 1974. vii+82 pp

\smallskip

\hskip -1 truecm 15. \  Y. Katznelson, D. Ornstein, {\it The differentiability of the conjugation 
of certain diffeomorphisms of the circle.} Erg. Th. and Dyn. Syst. {\bf 9},
643-680 (1989).

\smallskip 
\hskip -1 truecm 16. \ I. Liousse, {\it  Nombre de rotation dans les groupes de
Thompson généralisés, applications.} http ://hal.ccsd.cnrs.fr/ccsd-00004554.
\smallskip

\smallskip


\hskip -1 truecm 17. \  G. Margulis, {\it Free subgroups of the homeomorphism
 group of the circle.}  
 C. R. Acad. Sci. Paris Sér. I Math.  {\bf 331}(9), 669--674 (2000).

\smallskip


\hskip -1 truecm 18. \ H.  Minakawa, {\it  Classification of exotic circles of
  ${\rm   PL}\sb +(S\sp 1)$}.  Hokkaido Math. J. {\bf 26} (3), 685--697 (1997).

\smallskip


\hskip -1 truecm 19. \ H. Poincar\'e, {\it Oeuvres compl\`etes.}
 {\bf t.1}, 137-158.

\smallskip


\hskip -1 truecm 20. \   G. Rhin, {\it Approximants
  de Padé et mesures effectives d'irrationnalité}, séminiare de théorie des
nombres  Paris 1985-1986, Prog. math  {\bf 71}, 155-164 (1987).

\smallskip

\hskip -1 truecm 21.\  Y. Rubin, {\it Locally moving groups and reconstruction
problems}, Math. Appl. {\bf 354}, Kluwer Acad. Publ., Dordrecht (1996). 

\smallskip

\hskip -1 truecm 22.\  M. Stein,  {\it Groups of piecewise linear
  homeomorphisms}, Trans. A.M.S. {\bf 332},  477-514  (1992).

\vfill

\noindent Laboratoire Paul PAINLEV\'E,  U.M.R  8524  au CNRS 
\hfill\break
U.F.R. de Math\'ematiques, Universit\'e Lille I
\hfill\break
59655 Villeneuve d'Ascq Cédex, France

\end{document}